\documentclass[12pt]{amsart} 
\usepackage[mathscr]{eucal}
\usepackage{amsmath,amsfonts}

\parskip=\smallskipamount

\newtheorem{theorem}{Theorem}[section]

\newtheorem{lemma}[theorem]{Lemma}

\makeatletter
\@addtoreset{equation}{section}
\makeatother

\newcommand{\CC}{{\mathbb C}}
\newcommand{\NN}{{\mathbb N}}

\newcommand{\ZZ}{{\mathbb Z}}
\newcommand{\DD}{{\mathbb D}}

\newcommand{\FF}{{\mathbb F}}

\newcommand{\cA}{{\mathcal A}}

\newcommand{\cD}{{\mathcal D}}

\newcommand{\cF}{{\mathcal F}}

\newcommand{\cH}{{\mathcal H}}
\newcommand{\cK}{{\mathcal K}}
\newcommand{\cL}{{\mathcal L}}

\newcommand{\cP}{{\mathcal P}}
\newcommand{\cR}{{\mathcal R}}

\newcommand{\cT}{{\mathcal T}}
\newcommand{\fI}{{\mathfrak I}}
\newcommand{\cW}{{\mathcal W}}

\begin{document}
\pagestyle{plain}

\bigskip

\title{Tensor algebras and displacement structure.  \\
IV. Invariant kernels} 
\author{T. Banks} \author{T. Constantinescu} \author{Nermine El-Sissi} 

\address{Department of Mathematics \\
  University of Texas at Dallas \\
  Richardson, TX 75083}
\email{\tt banks@utdallas.edu} 
\address{Department of Mathematics \\
  University of Texas at Dallas \\
  Richardson, TX 75083} 
\email{\tt tiberiu@utdallas.edu}
\address{Department of Mathematics \\
  University of Texas at Dallas \\
  Richardson, TX 75083}
\email{\tt }

\noindent
\begin{abstract}
In this paper we investigate the class of 
invariant positive
definite kernels on the free semigroup on $N$ generators.
We provide a combinatorial description of the positivity of the kernel in 
terms of Dyck paths and then we find a displacement equation that encodes 
the invariance property of the kernel.

\end{abstract}

\maketitle

\section{Introduction}
In the previous parts of this paper, \cite{BC}, \cite{CJ}, 
there were considered algebraic and asymptotic properties of orthogonal
polynomials in several variables associated with a certain 
class of
positive definite kernels on the free semigroup on 
$N$ generators. These kernels were naturally associated with
the Cuntz-Toeplitz defining relations, $X_iX_j=\delta _{i,j}1$, 
$i,j=1,\ldots ,N$, but they are quite sparse (a lot
of zero entries), which makes their structure to be quite simple 
(see \cite{CJ} for more details).

In this note we consider a more general class of 
positive definite kernels, which are invariant under the action of the 
free semigroup on itself by concatenation, 
and our main goal is to find combinatorial descriptions of the positive
definiteness and of the invariance property.

In answering the first question we establish a connection with the 
combinatorics of Dyck paths. The invariance is then encoded into a 
displacement equation, and this allows the use of the tools of the 
displacement structure theory.

The paper is organized as follows. In Section ~2 we review some material
on orthogonal polynomials and introduce the moment kernel of a 
q-positive functional on the algebra of polynomials in several noncommuting
variables. Then we describe the main result about positive definite kernels
and Dyck paths in Theorem ~2.2. Two simple applications are given to the 
counting of paths in marine seismology and to the structure of the Markov
product introduced in \cite{Boz}.
Finally we discuss the connection between orthogonal polynomials
and displacement equations in our setting. The main result here is 
given by Theorem ~2.4. In Section ~3 we introduce several examples of 
positive definite invariant kernels. First we deal with a kernel 
involved in the dilation theory for arbitrary families of contractions and 
calculate the orthogonal polynomials associated with this kernel in 
Theorem ~3.1. Then we show in Theorem ~3.3 that the invariant kernels
are precisely the moment kernels associated with q-positive functionals 
on the algebra of polynomials in $N$ noncommuting 
isometric variables. Finally, we show that free products give many 
examples of invariant kernels. In Section ~4 we prove the main result about 
the displacement equation satisfied by an invariant kernel.

\section{Preliminaries}
Here we describe our setting, so that the paper can be 
read independently of \cite{BC} and \cite{CJ}.
We review material on orthogonal polynomials, positive definite kernels, and 
displacement structure. We also establish the description of 
positive definite kernels in terms of the combinatorics of Dyck paths.

\subsection{Moment kernels and orthogonal polynomials}
We introduce a class of positive definite kernels associated
with some linear functionals on algebras of polynomials.
Let $\cP _N$ be the algebra of polynomials in $N$ noncommuting
variables  $X_1,\ldots ,X_N$ with complex coefficients.
Each element $P\in \cP _N$ can be uniquely written 
in the form $P=\sum _{\sigma \in \FF _N^+}c_{\sigma }X_{\sigma }$,
with $c_{\sigma }\ne 0$ for finitely many $\sigma $'s, where 
$\FF _N^+$ denotes the unital free semigroup 
on $N$ generators $1,\ldots ,N$ and with lexicographic order $\preceq $.
Also $X_{\sigma }=X_{i_1}\ldots X_{i_k}$ for $\sigma =i_1\ldots i_k
\in \FF _N^+$. Instead of $\FF _1^+$ we use the standard notation
$\NN _0$.

We can view $\cP _N$ as the free product of $N$ copies of $\cP _1$:
$$\cP _N=\underbrace{\cP _1\star \ldots \star \cP _1}_{\mbox{$N$ times}}=
\CC\oplus 
\left(\oplus _{n\geq 1}
\oplus _{i_1\ne i_2,\ldots ,i_{n-1}\ne i_n}\cP ^0_{i_1}
\otimes \ldots \otimes \cP ^0_{i_n}\right),$$
where $\cP ^0_{i}$ is the set of polynomials in the variable $X_i$, 
$i=1,\ldots ,N$, without constant term.
We also notice 
that $\cP _N$ is isomorphic with the tensor algebra 
over $\CC ^N$, which is defined by the algebraic direct sum 
$$\cT (\CC ^N)=\oplus _{k\geq 0}(\CC ^N)^{\otimes k},$$
where  $(\CC ^N)^{\otimes k}$ denotes the $k$-fold tensor
product of $\CC ^N$ with itself.
If $\{e_1,\ldots ,e_N\}$ is the standard basis of
$\CC ^N$, then the set
$$\{1\}\cup \{e_{i_1}\otimes \ldots \otimes e_{i_k}  
\mid 1\leq i_1, \ldots ,i_k\leq N, k\geq 1\}
$$ 
is a  
basis of $\cT (\CC ^N)$. For $\sigma =i_1\ldots i_k$ we write
$e_{\sigma }$ instead of $e_{i_1}\otimes \ldots \otimes e_{i_k}$, 
and the mapping $X_{\sigma }\rightarrow e_{\sigma }$, $\sigma \in \FF _N^+$, 
extends to an isomorphism from  $\cP _N$ to $\cT (\CC ^N)$, 
hence $\cP _N \simeq \cT (\CC ^N)$.

There is a natural involution on 
$\cP _{2N}$ introduced as follows:
$$X^+_k=X_{N+k},\quad k=1,\ldots ,N,$$
$$X^+_l=X_{l-N},\quad l=N+1,\ldots ,2N;$$
on monomials, 
$$(X_{i_1}\ldots X_{i_k})^+=X^+_{i_k}\ldots X^+_{i_1},$$
and finally, if $Q=\sum _{\sigma \in \FF _{2N}^+}c_{\sigma }X_{\sigma }$,
then $Q^+=
\sum _{\sigma \in \FF _{2N}^+}\overline{c}_{\sigma }
X^+_{\sigma }$.

We say that $\cA \subset \cP _{2N}$
is  {\it symmetric} with respect to this involution if $P\in \cA $ implies 
$cP^+\in \cA $ for some $c\in \CC -\{0\}$. 
Then the quotient of $\cP _{2N}$
by the two-sided ideal generated by 
$\cA $ is an associative algebra $\cR (\cA )$.
Letting $\pi =\pi _{\cA }:\cP _{2N}\rightarrow 
\cR (\cA )$ denote the quotient map then the formula 
\begin{equation}\label{invo}
\pi  (P)^+=\pi (P^+)
\end{equation}
gives a well-defined involution on $\cR (\cA )$.
Thus $\cR (\cA )$ is a unital $*$-algebra and $\cA$ is called the set
of {\it defining relations}. 
A linear functional $\phi $
on $\cR (\cA )$ is called {\it q-positive} if 
$\phi (\pi (P)^+\pi (P))\geq 0$
for all $P\in \cP _N$. 
The {\em index set} $G(\cA )\subset \FF ^+_N$ of $\cA$ 
is  
chosen as follows:
if $\alpha \in G(\cA )$, choose the next
element in $G(\cA )$ to be the least $\beta \in \FF _N^+$ 
with the property that the 
elements $\pi (X_{\alpha '})$, $\alpha '\preceq \alpha $, 
and $\pi (X_{\beta })$ are linearly independent.
We will avoid the degenerate situation in which $\pi (1)=0$;
if we do so, then $\emptyset \in G(\cA )$. 
Define $F_{\alpha }=\pi (X_{\alpha })$ for $\alpha \in G(\cA )$.
The {\em moments} of $\phi $ are the complex numbers
\begin{equation}\label{mome}
s_{\alpha ,\beta }=
\phi (F^+_{\alpha }F_{\beta }), 
\quad \alpha ,\beta \in G(\cA ),
\end{equation}
and the {\em moment kernel } is defined by 
 $K_{\phi }(\alpha ,\beta )=
s_{\alpha ,\beta }$, $\alpha ,\beta \in G(\cA )$. 
Since $\phi $ is q-positive  on 
$\cR (\cA )$, $K_{\phi }$ 
is a positive definite kernel on $G(\cA )$.
However, $K_{\phi }$ does not determine $\phi $ uniquely. 

In \cite{BC} and \cite{CJ} the  focus was on moment kernels
associated with q-positive functionals on $\cR (\cA ^N_{CT})$, where 
$\cA ^N_{CT}=\{1-X^+_kX_k\mid k=1,\ldots ,N\}\cup
\{X^+_kX_l, k,l=1,\ldots ,N, k\ne l\}$. The relations 
$X^+_kX_l=0$, $k,l=1,\ldots ,N$, make the moment kernel to be sparse.
In this paper we analize the moment kernels of q-positive
functionals on $\cR (\cA ^N_O)$, where 
$\cA ^N_O=\{1-X^+_kX_k\mid k=1,\ldots ,N\}$. We have that 
$G(\cA ^N_O)=\FF ^+_N $.

The  
orthonormal polynomials 
associated with a strictly q-positive functional on $\cR (\cA )$
(that is, $\phi (\pi (P^+)\pi (P))>0$ for $\pi (P)\ne 0$) are 
introduced by the 
Gram-Schmidt procedure applied to the family 
$\{\pi (X_{\alpha })\}_{\alpha \in G(\cA )}$ of linearly independent
elements in the Hilbert space $\cH _{\phi }$ associated with $\phi $
by the Gelfand-Naimark-Segal construction. Thus, the 
orthonormal polynomials are
\begin{equation}\label{bond1}
\varphi _{\alpha  }=
\sum _{\beta \preceq \alpha }a_{\alpha ,\beta }\pi (X_{\beta }),
\quad a_{\alpha ,\alpha }>0.
\end{equation} 
The 
polynomials $\varphi _{\alpha  }$ are uniquely determined by the condition 
$a_{\alpha ,\alpha }>0$ and the orthonormality property
\begin{equation}\label{bond2}
\phi (\varphi ^+_{\beta  } \varphi _{\alpha  })=
\delta _{\alpha ,\beta },
\quad \alpha ,\beta \in G(\cA ).
\end{equation}

\subsection{Positive definite kernels and Dyck paths}
We will use several times a certain structure  (and parametrization)
of positive definite kernels on $\NN _0$.
For sake of completeness we briefly describe this structure here, 
while the details can be found in \cite{Co}. Also we discuss the 
connection with the combinatorics of Dyck paths.

For a contraction $\gamma \in \cL (\cH ,\cH ')$, that is, 
a linear bounded operator between the Hilbert spaces $\cH$ and 
$\cH '$ with $\|\gamma \|\leq 1$, we define the {\it defect operator}
$d_{\gamma }=(I-\gamma ^*\gamma )^{1/2}$, and the corresponding 
{\it defect space} $\cD _{\gamma }$, the closure of the range of
$d_{\gamma }$.  The {\it Julia operator} associated with $\gamma $
is defined by 
$$J(\gamma )=\left[\begin{array}{cc}
\gamma & d_{\gamma ^*} \\
d_{\gamma } & -\gamma ^*
\end{array}
\right];
$$
the Julia operator is unitary from $\cH \oplus \cD _{\gamma ^*}$
onto $\cH '\oplus \cD _{\gamma }$. This construction can be extended to
certain families of contractions as follows. Let 
$\{\gamma _{k,j}\}_{0\leq k\leq j}$ be a family of contractions satisfying 
the compatibility conditions: $(i)$ $\gamma _{k,k}=0$ for all $k\geq 0$ and 
$(ii)$ $\gamma _{k,j}\in 
\cL (\cD _{\gamma _{k+1,j}},\cD _{\gamma ^*_{k,j-1}})$.
Then the unitary operators $U_{k,j}$ are recursively defined
by: $U_{k,k}=I_1$ and for $k<j$,
$$\begin{array}{rcl}
U_{k,j}& =& \left(J(\gamma _{k,k+1})\oplus I_{j-k-1}\right)
\left(I_1\oplus J(\gamma _{k,k+2})\oplus I_{j-k-2}\right)\ldots
\left( I_{j-k-1}\oplus J(\gamma _{k,j})\right) \\
 & & \\
 & & \times  \left(U_{k+1,j}\oplus I_1\right);
\end{array}
$$
each $U_{k,j}$ is a $(j-k+1)\times (j-k+1)$
block matrix and $I_l$ denotes the identity $l\times l$ block matrix.

\begin{figure}[h]
\setlength{\unitlength}{3447sp}%
\begingroup\makeatletter\ifx\SetFigFont\undefined%
\gdef\SetFigFont#1#2#3#4#5{%
  \reset@font\fontsize{#1}{#2pt}%
  \fontfamily{#3}\fontseries{#4}\fontshape{#5}%
  \selectfont}%
\fi\endgroup%
\begin{picture}(5028,2364)(187,-1615)
{ \thinlines

\put(301,164){\circle{212}}
\put(2101,164){\circle{212}}
\put(5101,164){\circle{212}}

\put(601,-361){\framebox(600,600){}}
\put(2401,-361){\framebox(600,600){}}
\put(4201,-361){\framebox(600,600){}}
\put(1501,-811){\framebox(600,600){}}
\put(3301,-811){\framebox(600,600){}}
\put(2401,-1261){\framebox(600,600){}}

\put(301,164){\line( 1, 0){4800}}
\put(301,-286){\line( 1, 0){4800}}
\put(301,-736){\line( 1, 0){4800}}
\put(301,-1186){\line( 1, 0){4800}}
\put(676,-286){\line( 1, 1){450}}
\put(676,164){\line( 1,-1){450}}
\put(2476,164){\line( 1,-1){450}}
\put(4276,164){\line( 1,-1){450}}
\put(2476,-286){\line( 1, 1){450}}
\put(4276,-286){\line( 1, 1){450}}
\put(1576,-736){\line( 1, 1){450}}
\put(3376,-736){\line( 1, 1){450}}
\put(2476,-1186){\line( 1, 1){450}}
\put(1576,-286){\line( 1,-1){450}}
\put(3376,-286){\line( 1,-1){450}}
\put(2476,-736){\line( 1,-1){450}}
\put(301, 89){\line( 0, 1){450}}
\put(2101, 89){\line( 0, 1){450}}
\put(5101, 89){\line( 0, 1){450}}
\put(226,164){\line( 1, 0){ 75}}
\put(5101,164){\line( 1, 0){ 75}}
\put(301,539){\vector( 0,-1){300}}
\put(2101,539){\vector( 0,-1){300}}
\put(5101,239){\vector( 0, 1){300}}
}%
\put(676,389){$J(\gamma _{23})$}%

\put(2476,389){$J(\gamma _{12})$}%

\put(4276,389){$J(\gamma _{01})$}%

\put(1576,-61){$J(\gamma _{13})$}%

\put(3376,-61){$J(\gamma _{02})$}%

\put(2476,-1561){$J(\gamma _{03})$}%

\put(226,614){$A$}%

\put(2026,614){$C$}%

\put(5026,614){$B$}%

\end{picture}

\caption{\mbox{ Transmission line for $K(0,3)$ }}
\end{figure}

\begin{theorem}\label{structure}
Let $K$ be a positive definite kernel on the set $\NN _0$ with values 
in a Hilbert space $\cH $. Then there is a uniquely determined family of 
contractions satisfying the compatibility conditions $(i)$ and $(ii)$, 
and such that 
\begin{equation}\label{prima}
K(l,m)=K(l,l)^{1/2}\left(P_{\cH }U_{l,m}/\cH\right)K(m,m)^{1/2}, 
\quad l\leq m,
\end{equation}
where $P_{\cH }$ denotes the orthogonal projection on the space $\cH $.
\end{theorem}

For a proof see \cite{Co}. 
We shall say that $\{\gamma _{k,j}\}$ is the {\it family of parameters}
associated with the kernel $K$. Ocassionally we write 
$\gamma _{k,j}(K)$ in order to underline the dependence on $K$.
It is very useful to realize the above formula by a
so-called time varying transmission line; for 
$K(0,3)$ this is illustrated in Figure~1 (for simplicity, assume
$K(l,l)=I$ for all $l$).
Thus, if the identity operator $I$ is the input at $A$, then at $B$
we read off the expression of $K(0,3)$ in terms of the 
parameters $\gamma _{01}$, $\gamma _{02}$,
$\gamma _{03}$, $\gamma _{12}$, $\gamma _{13}$, $\gamma _{23}$
and their defects. Likewise, if the input at $C$ is the identity 
operator, then the output at $B$ is now the expression of $K(0,2)$
(for more details see \cite{Co}).

It was noticed in \cite{BC} that there is a simple connection between 
transmission lines as in Figure ~1 and Dyck (or Catalan)
paths.  We assume, again for simplicity, that 
$K(l,m)\in \CC$ for all $l,m$ and that $K(l,l)=1$ for all $l$. 
A Dyck path of length $2k$ is a path in the 
positive quadrant of the lattice $\ZZ ^2$ which starts at $(0,0)$, 
ends at $(2k,0)$, and consists of rise steps $\nearrow $
and fall steps $\searrow $ (see Figure ~2). For more 
information on Dyck paths and their combinatorics, see \cite{St}.  

\begin{figure}[h]
\setlength{\unitlength}{3447sp}%
\begingroup\makeatletter\ifx\SetFigFont\undefined%
\gdef\SetFigFont#1#2#3#4#5{%
  \reset@font\fontsize{#1}{#2pt}%
  \fontfamily{#3}\fontseries{#4}\fontshape{#5}%
  \selectfont}%
\fi\endgroup%
\begin{picture}(2424,1224)(289,-673)
{ \thinlines
\multiput(301,239)(8.98876,0.00000){268}{\makebox(1.6667,11.6667){\SetFigFont{5}{6}{\rmdefault}{\mddefault}{\updefault}.}}

\multiput(301,-61)(8.98876,0.00000){268}{\makebox(1.6667,11.6667){\SetFigFont{5}{6}{\rmdefault}{\mddefault}{\updefault}.}}

\multiput(301,-361)(8.98876,0.00000){268}{\makebox(1.6667,11.6667){\SetFigFont{5}{6}{\rmdefault}{\mddefault}{\updefault}.}}

\multiput(601,-661)(0.00000,9.02256){134}{\makebox(1.6667,11.6667){\SetFigFont{5}{6}{\rmdefault}{\mddefault}{\updefault}.}}

\multiput(901,-661)(0.00000,9.02256){134}{\makebox(1.6667,11.6667){\SetFigFont{5}{6}{\rmdefault}{\mddefault}{\updefault}.}}

\multiput(1201,-661)(0.00000,9.02256){134}{\makebox(1.6667,11.6667){\SetFigFont{5}{6}{\rmdefault}{\mddefault}{\updefault}.}}

\multiput(1501,-661)(0.00000,9.02256){134}{\makebox(1.6667,11.6667){\SetFigFont{5}{6}{\rmdefault}{\mddefault}{\updefault}.}}

\multiput(1801,-661)(0.00000,9.02256){134}{\makebox(1.6667,11.6667){\SetFigFont{5}{6}{\rmdefault}{\mddefault}{\updefault}.}}

\multiput(2101,-661)(0.00000,9.02256){134}{\makebox(1.6667,11.6667){\SetFigFont{5}{6}{\rmdefault}{\mddefault}{\updefault}.}}

\multiput(2401,-661)(0.00000,9.02256){134}{\makebox(1.6667,11.6667){\SetFigFont{5}{6}{\rmdefault}{\mddefault}{\updefault}.}}

\put(301,-661){\framebox(2400,1200){}}

\put(301,-661){\line( 1, 1){600}}
\put(901,-61){\line( 1,-1){300}}
\put(1201,-361){\line( 1, 1){600}}
\put(1801,239){\line( 1,-1){900}}
}%

\end{picture}

\caption{\mbox{ A Dyck path of length $8$ }}
\end{figure}

Let $\cD _k$ 
be the set of Dyck paths of length $2k$ and 
let $\cA _k$ be the set of points $(l,q)$, $q>0$, 
with the property that there exists ${\bf p}\in \cD _k$ with 
$(l,q)\in {\bf p}$. It is seen that 
$$\cA _k=\left\{(j+i,j-i)\mid 0\leq i<j\leq k\right\}.$$
Also, we notice that if  ${\bf p}\in \cD _k$ and 
$x=(l,q)\in {\bf p}$, then there are only four types of behaviour of $
{\bf p}$
about $x$: (I) a rise step followed by a fall step; 
(II) a fall step followed by a rise step; (III)
two consecutive rise steps; (IV) two consecutive fall steps (see
Figure ~3). 

\begin{figure}[h]
\setlength{\unitlength}{3447sp}%
\begingroup\makeatletter\ifx\SetFigFont\undefined%
\gdef\SetFigFont#1#2#3#4#5{%
  \reset@font\fontsize{#1}{#2pt}%
  \fontfamily{#3}\fontseries{#4}\fontshape{#5}%
  \selectfont}%
\fi\endgroup%
\begin{picture}(4224,966)(589,-1015)
{ \thinlines
{\put(601,-661){\line( 1, 1){300}}
}%
{\put(901,-361){\line( 1,-1){300}}
}%
{\put(1801,-361){\line( 1,-1){300}}
}%
{\put(2101,-661){\line( 1, 1){300}}
}%
{\put(3001,-661){\line( 1, 1){600}}
}%
{\put(4201,-61){\line( 1,-1){600}}
}%
\put(826,-961){\makebox(0,0)[lb]{\smash{{\SetFigFont{12}{14.4}{\rmdefault}{\mddefault}{\updefault}{$(I)$}%
}}}}
\put(2026,-961){\makebox(0,0)[lb]{\smash{{\SetFigFont{12}{14.4}{\rmdefault}{\mddefault}{\updefault}{$(II)$}%
}}}}
\put(3226,-961){\makebox(0,0)[lb]{\smash{{\SetFigFont{12}{14.4}{\rmdefault}{\mddefault}{\updefault}{$(III)$}%
}}}}
\put(4276,-961){\makebox(0,0)[lb]{\smash{{\SetFigFont{12}{14.4}{\rmdefault}{\mddefault}{\updefault}{$(IV)$}%
}}}}
\put(850, -600){$x$}
\put(2050, -500){$x$}
\put(3320, -500){$x$}
\put(4300, -500){$x$}

}
\end{picture}

\caption{\mbox{ Behaviour of a Dyck path about a vertex $x\in \cA _k$}}
\end{figure}

Consequently, 
for each pair $i,j$ with $0\leq i<j\leq k$ we define the function
$a_{i,j}:\cD _k\rightarrow \CC $, 
$$a_{i,j}({\bf p})=\left\{\begin{array}{ccl}
1 &\mbox{if} & x=(j+i,j-i)\notin {\bf p}; \\
\gamma _{i,j} &\mbox{if} & x=(j+i,j-i)\in {\bf p} \quad \mbox{and (I) holds};\\
-\gamma ^*_{i,j} &\mbox{if} & x=(j+i,j-i)\in {\bf p} \quad \mbox{and (II) holds};\\
d_{i,j} &\mbox{if} & x=(j+i,j-i)\in {\bf p} \quad \mbox{and
either (III) or (IV) holds}.
\end{array}
\right.
$$
Let ${\bf p}$ be a Dyck path in $\cD _k$ such that $(2l,0)\in {\bf p}$.
The restriction of ${\bf p}$ from $(2l,0)$ to $(2k,0)$ is called
a {\it Dyck subpath starting at }$(2l,0)$ in $\cD _k$  and denote by  
$\cD ^l_k$ the set of all these subpaths. There is a bijection 
between $\cD ^l_k$ and $\cD _{k-l}$ so that the number of elements
in $\cD ^l_k$ is given by the Catalan number 
$C_{k-l}=\displaystyle\frac{1}{k-l+1}\left(
\begin{array}{c}
2(k-l) \\
k-l
\end{array}
\right);
$
also, $\cD ^0_k=\cD _k$.
If ${\bf q}\in \cD ^l_k$ then there could be many Dick paths 
whose restrictions at $(2l,0)$ coincide with ${\bf q}$. However, we 
notice that if ${\bf p}_1$ and ${\bf p}_2$ are two such Dych paths, then 
$a_{i,j}({\bf p}_1)=a_{i,j}({\bf p}_2)$ for $j+i>2l$. We will write 
$a_{i,j}({\bf q})$ in order to denote this common value. 

Now we can rewrite \eqref{prima} as a cumulant type formula. 
In fact, we can establish a 
certain connection with free cumulants (see \cite{VDN}, \cite{Sp}), 
which will be explored elsewhere.
\begin{theorem}\label{suma}
Let $K$ be a positive definite kernel on the set $\NN _0$ with scalar
values and $K(l,l)=1$ for all $l$.
Then, for $l<m$,
\begin{equation}\label{fsuma}
K(l,m)=\sum _{{\bf q}\in \cD ^l_m}\prod _{l\leq i<j\leq m}a_{i,j}({\bf q}).
\end{equation}
\end{theorem}

\begin{proof}
Formula \eqref{fsuma} is a direct consequence of Theorem ~\ref{structure} and 
the straightforward way in which we identify paths in a transmission line with 
Dyck paths.
\end{proof}

Formula \eqref{fsuma} looks quite intriguing. 
There is a well-established connection between continued fraction
expansions and combinatorics of Dyck paths, see for 
instance \cite{La},  still \eqref{fsuma}
comes from the only requirement that the kernel is positive.
A brief application of this result concerns the counting of 
paths in marine seismology. One has a layered medium with a perfect 
reflection at the $0$-interface (see Figure ~4).

\begin{figure}[h]
\setlength{\unitlength}{3947sp}%
\begingroup\makeatletter\ifx\SetFigFont\undefined%
\gdef\SetFigFont#1#2#3#4#5{%
  \reset@font\fontsize{#1}{#2pt}%
  \fontfamily{#3}\fontseries{#4}\fontshape{#5}%
  \selectfont}%
\fi\endgroup%
\begin{picture}(3742,1497)(589,-973)
\thinlines
{\put(601,239){\line( 1, 0){2400}}
}%
{\put(601,-361){\line( 1, 0){2400}}
}%
{\put(601,-961){\line( 1, 0){2400}}
}%
{\put(901,239){\line( 1,-2){300}}
\put(1201,-361){\line( 1, 2){300}}
\put(1501,239){\line( 1,-2){300}}
\put(1801,-361){\vector( 1, 2){300}}
}%
\put(3301,239){\makebox(0,0)[lb]{\smash{{\SetFigFont{12}{14.4}{\rmdefault}{\mddefault}{\updefault}{\mbox{$0$-interface}}%
}}}}
\put(3301,-286){\makebox(0,0)[lb]{\smash{{\SetFigFont{12}{14.4}{\rmdefault}{\mddefault}{\updefault}{\mbox{$1$st interface}}%
}}}}
\put(3301,-886){\makebox(0,0)[lb]{\smash{{\SetFigFont{12}{14.4}{\rmdefault}{\mddefault}{\updefault}{\mbox{$2$nd interface}}%
}}}}
\put(901,389){\makebox(0,0)[lb]{\smash{{\SetFigFont{12}{14.4}{\rmdefault}{\mddefault}{\updefault}{$A$}%
}}}}
\end{picture}%

\caption{\mbox{ A trajectory through a layered medium in $4$ units of time}}
\end{figure}

A unit impulse strikes at $A$, at time zero, and it propagates downwards 
through the medium. At each interface, the impulse is partially
reflected and partially transmitted to the next layer. It is a consequence of 
Theorem ~\ref{suma} that the number of possible paths the impulse can take 
in order to return back to the $0$-interface in $2n$ units of time, 
is precisely given by the Catalan number $C_n$.

\bigskip
We conclude this subsection with an application of the 
transmission line interpretation of 
Theorem ~\ref{structure} which provides a simple, conceptual proof of 
a result in \cite{Boz}. Thus, let $A_1$, $A_2$ be two sets such 
that $A_1\cap A_2=\{a\}$ and $K_1$, $K_2$ be positive definite
kernels on $A_1$, respectively $A_2$, such that 
$K_1(a,a)=K_2(a,a)$. The Markov product
of the kernels $K_1$ and $K_2$ is a hermitian kernel
$K$ on $A_1 \cup A_2$ defined in \cite{Boz} by the rules:
$$\begin{array}{l}
(1) \quad K\mid_{A_j\times A_j}=K_j, \quad j=1,2;\\
 \\
(2) \quad K(a_1,a_2)=K_1(a_1,a)K_2(a,a_2), \quad a_1\in A_1, a_2\in A_2;\\
 \\
(3) \quad K(a_1,a_2)=K(a_2,a_1)^*.
\end{array}
$$
For our purpose we can restrict to the case of finite sets $A_1$, $A_2$, 
$A_1=\{0,\ldots ,n\}$, $A_2=\{-m,\ldots ,0\}$. Let 
$\{\gamma _{k,j}(K_1)\}_{-m\leq k\leq j\leq 0}$
and $\{\gamma _{k,j}(K_2)\}_{0\leq k\leq j\leq n}$
be the parameters associated with $K_1$, respectively $K_2$.
The fact that the Markov product is positive definite was proved in 
\cite{Boz}. In addition, we provide here the structure of its 
associated parameters.

\begin{theorem}\label{boz}
The Markov product of two positive definite 
kernels $K_1$ and $K_2$ is a positive definite kernel with parameters
$\{\gamma _{k,j}\}_{-m\leq k\leq j\leq n}$ given by:
$$\gamma _{k,j}=\left\{
\begin{array}{ccl}
\gamma _{k,j}(K_1) &\mbox{if} & -m\leq k\leq j\leq 0; \\
\gamma _{k,j}(K_2) &\mbox{if} & 0\leq k\leq j\leq n; \\
0 & & \mbox{otherwise}.
\end{array}\right.
$$
\end{theorem}
\begin{proof}
The transmission line of $K(a_1,a_2)$, 
$a_1\in A_1$, $a_2\in A_2$, looks like in Figure ~5.

\begin{figure}[h]
\setlength{\unitlength}{2447sp}%
\begingroup\makeatletter\ifx\SetFigFont\undefined%
\gdef\SetFigFont#1#2#3#4#5{%
  \reset@font\fontsize{#1}{#2pt}%
  \fontfamily{#3}\fontseries{#4}\fontshape{#5}%
  \selectfont}%
\fi\endgroup%
\begin{picture}(8592,2856)(259,-2173)
\thinlines

\put(601,-361){\framebox(600,600){}}
\put(1501,-811){\framebox(600,600){}}
\put(2401,-361){\framebox(600,600){}}
\put(3301,-811){\framebox(600,600){}}
\put(5101,-811){\framebox(600,600){}}
\put(4201,-361){\framebox(600,600){}}
\put(6001,-361){\framebox(600,600){}}
\put(7801,-361){\framebox(600,600){}}
\put(6901,-811){\framebox(600,600){}}

\put(301,164){\line( 1, 0){8400}}
\put(3301,164){\vector( 1, 0){600}}
\put(301,464){\vector( 0,-1){300}}
\put(8701,164){\vector( 0, 1){300}}
\put(301,-286){\line( 1, 0){3075}}

\put(2401,-1261){\framebox(600,600){}}
\put(4201,-1261){\framebox(600,600){}}
\put(6001,-1261){\framebox(600,600){}}
\put(3301,-1711){\framebox(600,600){}}
\put(5101,-1711){\framebox(600,600){}}
\put(4201,-2161){\framebox(600,600){}}

\put(301,-736){\line( 1, 0){2175}}
\put(301,-1186){\line( 1, 0){2175}}

\put(2476,-736){\makebox(1.6667,11.6667){\SetFigFont{5}{6}{\rmdefault}{\mddefault}{\updefault}.}}

\put(2476,-736){\line( 1,-1){450}}
\put(2476,-1186){\line( 1, 1){450}}
\put(2926,-1186){\line( 1, 0){450}}
\put(3376,-1186){\line( 1,-1){450}}
\put(3376,-286){\line( 1,-1){450}}

\put(3826,-286){\makebox(1.6667,11.6667){\SetFigFont{5}{6}{\rmdefault}{\mddefault}{\updefault}.}}

\put(3451,-736){\line( 5, 6){375}}
\put(2926,-736){\line( 1, 0){525}}
\put(3826,-286){\line( 1, 0){4875}}
\put(3826,-736){\line( 1, 0){450}}
\put(3376,-1636){\line( 1, 1){450}}
\put(3826,-1186){\line( 1, 0){450}}
\put(4276,-736){\line( 1,-1){450}}
\put(4276,-1186){\line( 1, 1){450}}
\put(4726,-736){\line( 1, 0){3975}}
\put(4726,-1186){\line( 1, 0){450}}
\put(5176,-1186){\line( 1,-1){450}}
\put(3826,-1636){\line( 1, 0){450}}
\put(4276,-1636){\line( 1,-1){450}}
\put(4276,-2086){\line( 1, 1){450}}
\put(5251,-1636){\line( 5, 6){375}}
\put(4726,-1636){\line( 1, 0){525}}
\put(5626,-1186){\line( 1, 0){3075}}

\put(676,-286){\line( 1, 1){450}}
\put(2476,-286){\line( 1, 1){450}}
\put(4276,-286){\line( 1, 1){450}}
\put(6076,-286){\line( 1, 1){450}}
\put(7876,-286){\line( 1, 1){450}}
\put(1576,-736){\line( 1, 1){450}}
\put(5176,-736){\line( 1, 1){450}}
\put(6076,-1186){\line( 1, 1){450}}
\put(676,164){\line( 1,-1){450}}
\put(2476,164){\line( 1,-1){450}}
\put(1576,-286){\line( 1,-1){450}}
\put(4276,164){\line( 1,-1){450}}
\put(6076,164){\line( 1,-1){450}}
\put(7876,164){\line( 1,-1){450}}
\put(5176,-286){\line( 1,-1){450}}
\put(6976,-286){\line( 1,-1){450}}
\put(6076,-736){\line( 1,-1){450}}
\put(6976,-736){\line( 1, 1){450}}

\put(451,389){$A$}%
\put(8851,389){$B$}%
\put(3526,389){$C$}%
\end{picture}

\caption{\mbox{ Transmission line for Markov products}}
\end{figure}

The Julia operator of $0$ is $\left[\begin{array}{cc}
0 & I \\
I & 0
\end{array}
\right]$
and therefore the central block of Julia operators of $0$ acts like
a barrier. There is only one place for the signal to propagate from the 
left to right and that is the upmost wire. What comes through that wire 
is exactly $K_2(a,a_2)$. The transmission line to the right of $C$ will 
produce $K_1(a_1,a)$ and all together will get the product
$K_2(a,a_2)K_1(a_1,a)$.
\end{proof}

\subsection{Displacement structure and orthogonal polynomials}
The displacement structure of a family $\{R(t)\}$
of matrices is encoded by an equation of the form
$$R(t)-F(t)R(t+1)F(t)^*=G(t)J(t)G^*(t),$$
where $F(t)$, $G(t)$ are the so-called {\it generators}
and $J(t)$ is a signature matrix (usually, $J(t)=I_p\oplus -I_q$,
for some fixed $p,q$). The main feature in the use of displacement 
structure is that under suitable conditions on generators, the 
Gaussian elimination for $R(t)$ can be performed 
at the level of generators. This leads to faster algorithms for 
factorization of $R(t)$ and to useful
lattice structures associated with these matrices.

There is a remarkable connection 
between orthogonal polynomials on the unit circle
(when the moment kernel is Toeplitz) and displacement 
structure, as described in \cite{Ka}. For our purpose it is convenient
to obtain a similar connection in our more general setting. We discuss 
in details the following situation
(with the notation introduced in Subsection ~2.1):
$N=1$ and $\cA =\emptyset $, so that $\cR (\cA )=\cP _2$ 
and $G(\cA )=\NN _0$. 
The moment kernel of a q-positive 
functional $\phi $ on $\cP _2$ is $K(n,m)=\phi ((X_1^n)^+X_1^m)$, 
$n,m\in \NN _0 $, and there is no additional restriction on 
$K_{\phi }$ other then being positive definite. So, in a certain sense, 
this is the most general possible situation. 
Next assume $\phi $ is strictly q-positive (we say in this case that 
the moment kernel is strictly positive definite). It was showed in 
\cite{CJ} that the orthonormal polynomials
associated with $\phi $ obey the recurrence relation:

\begin{equation}\label{lazero}
\varphi _0(X_1,l)=\varphi _0^{\sharp }(X_1,l)=s_{l,l}^{-1/2}, 
\quad l\in \NN _0,
\end{equation}
and for $n\geq 1$, $l\in \NN _0$,
\begin{equation}\label{primarelatie}
\varphi _n(X_1,l)=\frac{1}{d_{l,n+l}}
\left( X_1\varphi _{n-1}(X_1,l+1)-
\gamma _{l,n+l}\varphi ^{\sharp }_{n-1}(X_1,l)\right),
\end{equation}
\begin{equation}\label{adouarelatie}
\varphi ^{\sharp}_n(X_1,l)=\frac{1}{d_{l,n+l}}
\left(-\overline{\gamma }_{l,n+l}X_1\varphi _{n-1}(X_1,l+1)+
\varphi ^{\sharp }_{n-1}(X_1,l)\right),
\end{equation}
where $\varphi _n(X_1)=\varphi _n(X_1,0)$ and $\{\gamma _{k,j}\}$
is the family of parameters associated with the moment kernel
$K_{\phi }$.

We now describe the displacement structure of the kernel $K_{\phi }$.
For each $n\geq 0$ we introduce the following elements (the 
generators of the relevant displacement equations): 
the $(n+1)\times (n+1)$ matrix
$$F_n(t)=\left[
\begin{array}{cccccc}
0 & & & &  \\
1 & 0 &  & \mbox{\bf \Large{0}} & \\
  & 1 & & &  \\
  & \mbox{\bf \Large{0}} & \ddots & \ddots &  \\
  &   & & 1 & 0
\end{array}
\right], \quad t\in \NN _0,
$$
and the $2\times 2$ matrix
$$J(t)=\left[
\begin{array}{cc}
1 & 0 \\
0 & -1
\end{array}
\right], \quad t\in \NN _0;
$$
then for $t\in \NN _0$, we introduce the $(n+1)\times 2$ matrix
$$G_n(t)=s^{-1/2}_{t,t}\left[
\begin{array}{cc}
s_{t,t} & 0 \\
\overline{s}_{t,t+1} & \overline{s}_{t,t+1} \\
\vdots & \vdots \\
\overline{s}_{t,t+n} & \overline{s}_{t,t+n}
\end{array}
\right].
$$
It was showed in \cite{CSK} that the displacement equation
\begin{equation}\label{displ}
R_n(t)-F_n(t)R_n(t+1)F_n(t)^*=G_n(t)J(t)G_n(t)^*, \quad t\in \NN _0 ,
\end{equation}
has a unique solution given by 
$R_n(t)=\left[s_{k,j}\right]_{t\leq k,j\leq t+n}$.
For this reason we say that the kernel $K_{\phi }$ has displacement 
structure.

If $K_{\phi }$ 
is a Toeplitz kernel then there is a strong connection between the 
displacement structure of its inverse and the orthogonal 
polynomials associated with $\phi $. In the general case there is 
a trade-off. Certainly, the more general formulae are somewhat 
obscured by the necessary use of additional indices. On the other hand, 
the general case reveals some features obscured by the additional 
symmetries of the Toeplitz case.

The orthonormal polynomials of $\phi $ have the expansion
$$\varphi _n(X_1,l)=\sum _{k=0}^na^l_{n,k}X^k_1,$$
with $\left(a^l_{n,n}\right)^{-1}=s^{1/2}_{l+n,l+n}
\prod _{k=1}^nd_{l+n-k,l+n}$, and similarly, the polynomials
$\varphi ^{\sharp }_n$ have the expansion
$$\varphi ^{\sharp }_n(X_1,l)=\sum _{k=0}^nb^l_{n,k}X^k_1,$$
with $\left(b^l_{n,0}\right)^{-1}=s^{1/2}_{l+n,l+n}
\prod _{k=1}^nd_{l,l+k}$.
It follows from the 
proof of Theorem ~3.2 in \cite{CJ} that 

\begin{equation}\label{alb1}
R_n(t)\left[
\begin{array}{c}
a^t_{n,0} \\
\vdots \\
a^t_{n,n-1} \\
a^t_{n,n}
\end{array}
\right]
=
\left[
\begin{array}{c}
0  \\
\vdots \\
0 \\
\left(a^t_{n,n}\right)^{-1}
\end{array}
\right]
\end{equation}
and 
\begin{equation}\label{alb2}
R_n(t)\left[
\begin{array}{c}
b^t_{n,0} \\
b^t_{n,1} \\
\vdots \\
b^t_{n,n}
\end{array}
\right]
=
\left[
\begin{array}{c}
\left(b^t_{n,n}\right)^{-1}\\
0  \\
\vdots \\
0 
\end{array}
\right].
\end{equation}

We then define 
$$H_n(t)=\left[
\begin{array}{cccc}
\overline{a}^{t+1}_{n,0} &
\ldots &
\overline{a}^{t+1}_{n,n-1} &
a^{t+1}_{n,n} \\
\overline{b}^{t}_{n,1} &
\ldots &
\overline{b}^{t}_{n,n} &
0
\end{array}
\right]
$$
and obtain the main result of this subsection.
\begin{theorem}\label{DSAOP}
The family $\{R_n(t)^{-1}\}_{t\in \NN _0 }$ is the solution of the 
displacement equation
$$R_n(t+1)^{-1}-F_n(t)^*R_n(t)^{-1}F_n(t)=H_n(t)^*J(t)H_n(t),\quad 
t\in \NN _0.$$
\end{theorem}
\begin{proof}

We define $K_n(t)=\left[
\begin{array}{cc}
0 & 0 \\
0 & -s^{1/2}_{t,t}b^t_{n,0}
\end{array}
\right]$, $t\in \NN _0$, and we have to show that 

\begin{equation}\label{albastru1}
F_n(t)R_n(t+1)H_n(t)^*+G_n(t)J(t)K_n(t)^*=0
\end{equation}
and 
\begin{equation}\label{albastru2}
H_n(t)R_n(t+1)H_n(t)^*+K_n(t)J(t)K_n(t)^*=J(t).
\end{equation}
From \eqref{alb2} we deduce that 
$$b^t_{n,0}\left[
\begin{array}{c}
\overline{s}_{t,t+1} \\
\vdots \\
\overline{s}_{t,t+n}
\end{array}\right]
+
\left[
\begin{array}{ccc}
s_{t+1,t+1} & \ldots & s_{t+1,t+n} \\
\vdots & \ddots & \\
\overline{s}_{t+1,t+n} & \ldots & s_{t+n,t+n}
\end{array}\right]
\left[
\begin{array}{c}
b^t_{n,1} \\
\vdots   \\
b^t_{n,n}
\end{array}
\right]=0
$$
which implies
$$R_n(t+1)\left[
\begin{array}{c}
b^t_{n,1} \\
\vdots \\
b^t_{n,n} \\
0
\end{array}
\right]
=
\left[
\begin{array}{c}
-b^t_{n,0}\overline{s}_{t,t+1} \\
\vdots \\
-b^t_{n,0}\overline{s}_{t,t+n} \\
{\mbox{\Large{*}}}
\end{array}
\right],
$$
where $\mbox{\Large{*}}$ denotes an entry whose actual value does not play
any role here. Therefore, using the previous relation and 
\eqref{alb1}, we deduce
$$R_n(t+1)H_n(t)^*=
\left[
\begin{array}{cc}
0 & -b^t_{n,0}\overline{s}_{t,t+1} \\
\vdots & \vdots \\
0 & -b^t_{n,0}\overline{s}_{t,t+n} \\
\left(a^{t+1}_{n,n}\right)^{-1} & \mbox{\Large{*}}
\end{array}
\right]
$$
and 
$$\begin{array}{rl}
F_n(t)R_n(t+1)H_n(t)^*+G_n(t)J(t)K_n(t)^* &  \\
 & \\
 & \!\!\!\!\!\!\!\!\!\!\!\!\!\!\!\!
 \!\!\!\!\!\!\!\!\!\!\!\!\!\!\!\!
 \!\!\!\!\!\!\!\!\!\!\!\!\!\!\!\!
 \!\!\!\!\!\!\!\!\!\!\!\!\!\!\!\!
 \!\!\!\!\!\!\!\!\!\!\!\!\!\!\!\!
 \!\!\!\!\!\!\!\!\!\!\!\!\!\!\!\!
 \!\!\!\!\!\!\!\!\!\!\!\!\!\!\!\!
=-b^t_{n,0}\left[
\begin{array}{cc}
0 & 0 \\
0 & \overline{s}_{t,t+1} \\
\vdots & \vdots \\
0 & \overline{s}_{t,t+n}
\end{array}
\right]+s^{-1/2}_{t,t}\left[
\begin{array}{cc}
s_{t,t} & 0 \\
\overline{s}_{t,t+1} & \overline{s}_{t,t+1} \\
\vdots & \vdots \\
\overline{s}_{t,t+n} & \overline{s}_{t,t+n} 
\end{array}
\right]
\left[
\begin{array}{cc}
1 & 0 \\
0 & -1
\end{array}
\right]
\left[
\begin{array}{cc}
0 & 0 \\
0 & -s^{1/2}_{t,t}b^t_{n,0}
\end{array}
\right] \\
 & \\
 & \!\!\!\!\!\!\!\!\!\!\!\!\!\!\!\!
 \!\!\!\!\!\!\!\!\!\!\!\!\!\!\!\!
 \!\!\!\!\!\!\!\!\!\!\!\!\!\!\!\!
 \!\!\!\!\!\!\!\!\!\!\!\!\!\!\!\!
 \!\!\!\!\!\!\!\!\!\!\!\!\!\!\!\!
 \!\!\!\!\!\!\!\!\!\!\!\!\!\!\!\!
 \!\!\!\!\!\!\!\!\!\!\!\!\!\!\!\!
=0.
\end{array}
$$
In order to obtain \eqref{albastru2} we calculate:
$$\begin{array}{rl}
H_n(t)R_n(t+1)H_n(t)^* & \\
 & \!\!\!\!\!\!\!\!\!\!\!\!\!\!\!\!
 \!\!\!\!\!\!\!\!\!\!\!\!\!\!\!\!
\!\!\!\!\!\!\!\!
=\left[
\begin{array}{cccc}
\overline{a}^{t+1}_{n,0} &
\ldots &
\overline{a}^{t+1}_{n,n-1} &
a^{t+1}_{n,n} \\
\overline{b}^{t}_{n,1} &
\ldots &
\overline{b}^{t}_{n,n} &
0
\end{array}
\right]
\left[
\begin{array}{cc}
0 & -b^t_{n,0}\overline{s}_{t,t+1} \\
\vdots & \vdots \\
0 & -b^t_{n,0}\overline{s}_{t,t+n} \\
\left(a^{t+1}_{n,n}\right)^{-1} & \mbox{\Large{*}}
\end{array}
\right] \\
 & \\
 & \!\!\!\!\!\!\!\!\!\!\!\!\!\!\!\!
\!\!\!\!\!\!\!\!\!\!\!\!\!\!\!\!
\!\!\!\!\!\!\!\!
=\left[
\begin{array}{cc}
1 & \mbox{\Large{*}} \\
0 & -b^t_{n,0}s_{t,t}(\overline{b}^t_{n,1}
\overline{s}_{t,t+1}+\ldots 
+\overline{b}^t_{n,n}
\overline{s}_{t,t+n})
\end{array}
\right].
\end{array}
$$
The fact that the matrix 
$H_n(t)R_n(t+1)H_n(t)^*$ is selfadjoint makes the north-east corner 
(the $\mbox{\Large{*}}$ entry) of the above matrix equal to $0$. Also, 
formula \eqref{alb1} implies that 
$$
\overline{b}^t_{n,1}
\overline{s}_{t,t+1}+\ldots 
+\overline{b}^t_{n,n}
\overline{s}_{t,t+n}
=\left(b^t_{n,0}\right)^{-1}-s_{t,t}b^t_{n,0}.$$
In conclusion, 
$$\begin{array}{rl}
H_n(t)R_n(t+1)H_n(t)^*+K_n(t)J(t)K_n(t)^* &  \\
 & \\
 & \!\!\!\!\!\!\!\!\!\!\!\!\!\!\!\!
\!\!\!\!\!\!\!\!\!\!\!\!\!\!\!\!
\!\!\!\!\!\!\!\!\!\!\!\!\!\!\!\!
\!\!\!\!\!\!\!\!\!\!\!\!\!\!\!\!
\!\!\!\!\!\!\!\!\!\!\!\!\!\!\!\!
=\left[
\begin{array}{cc}
1 & 0  \\
0 & -b^t_{n,0}\left(\left(b^t_{n,0}\right)^{-1}-s_{t,t}b^t_{n,0}\right)
\end{array}
\right]
+\left[
\begin{array}{cc}
0 & 0 \\
0 & s_{t,t}\left(b^t_{n,0}\right)^2
\end{array}
\right] \\
 & \\
 & \!\!\!\!\!\!\!\!\!\!\!\!\!\!\!\!
\!\!\!\!\!\!\!\!\!\!\!\!\!\!\!\!
\!\!\!\!\!\!\!\!\!\!\!\!\!\!\!\!
\!\!\!\!\!\!\!\!\!\!\!\!\!\!\!\!
\!\!\!\!\!\!\!\!\!\!\!\!\!\!\!\!
 =J(t).
\end{array}
$$
From \eqref{albastru1} and \eqref{albastru2}
we deduce that
$$\left[
\begin{array}{cc}
F_n(t) & G_n(t) \\
H_n(t) & K_n(t) 
\end{array}
\right]
\left[
\begin{array}{cc}
R_n(t+1) & 0 \\
 0 & J(t) 
\end{array}
\right]
\left[
\begin{array}{cc}
F_n(t) & G_n(t) \\
H_n(t) & K_n(t) 
\end{array}
\right]^*=
\left[
\begin{array}{cc}
R_n(t) & 0 \\
 0 & J(t) 
\end{array}
\right]
$$
and a Schur complement argument implies
$$\left[
\begin{array}{cc}
F_n(t) & G_n(t) \\
H_n(t) & K_n(t) 
\end{array}
\right]^*
\left[
\begin{array}{cc}
R_n(t)^{-1} & 0 \\
 0 & J(t) 
\end{array}
\right]
\left[
\begin{array}{cc}
F_n(t) & G_n(t) \\
H_n(t) & K_n(t) 
\end{array}
\right]=
\left[
\begin{array}{cc}
R_n(t+1)^{-1} & 0 \\
 0 & J(t) 
\end{array}
\right].
$$
In particular, we deduce that 
$$R_n(t+1)^{-1}-F_n(t)^*R_n(t)^{-1}F_n(t)=H_n(t)^*J(t)H_n(t),\quad 
t\in \NN _0.$$
\end{proof}

\section{Invariant kernels}
In this section we introduce invariant kernels and provide 
several examples. Let $\cH $ be a Hilbert space. A positive definite 
kernel $K:\FF ^+_N\times \FF ^+_N\rightarrow \cL (\cH )$ is 
called {\it invariant} (under the action of $\FF ^+_N$ on itself
by concatenation) if
\begin{equation}\label{invariance}
K(\tau \sigma ,\tau \sigma ')=K(\sigma ,\sigma '), \quad \tau, \sigma ,
\sigma '\in \FF ^+_N.
\end{equation} 

The invariant Kolmogorov decomposition theorem, \cite{Pa}, provides
a certain structure of an invariant kernel. Thus, we can define a Hilbert
space $\cK $, an operator $V\in \cL (\cH ,\cK )$, and an isometric
representation $U$ of $\FF ^+_N$ on $\cK $ such that
\begin{equation}\label{reprez}
K(\sigma ,\tau )=V^*U(\sigma )^*U(\tau )V, \quad \sigma ,
\tau \in \FF ^+_N,
\end{equation} 
and the set $\left\{U(\sigma Vh \mid \sigma \in \FF ^+_N, h\in \cH \right\}$
is total in $\cK $. An application of Theorem ~2.1 gives a family of 
contractions $\{\gamma _{\sigma ,\tau }\mid \sigma ,\tau \in \FF ^+_N, 
\sigma \preceq \tau \}$ satisfying the compatibility conditions:
$(i)$ $\gamma _{\sigma ,\sigma }=0$ for $\sigma \in \FF ^+_N$ 
and $(ii)$ $\gamma _{\sigma ,\tau }\in 
\cL (\cD _{\gamma _{\sigma +1,\tau }},\cD _{\gamma ^*_{\sigma ,\tau -1}})$,
where $\tau -1$ denotes the predecessor of $\tau $ with respect to the 
lexicographic order on $\FF ^+_N$ and $\sigma +1$ denotes the 
succesor of $\sigma $. Using Theorem ~1.6.1 in \cite{Co}, 
we can describe $V$ and $U$ in \eqref{reprez} in terms
of the parameters $\gamma _{\sigma ,\tau }$. However, at this stage
it is not clear how to translate the invariance property of $K$ into
 an invariance property of the parameters 
$\gamma _{\sigma ,\tau }$. In order to deal with this issue we first
discuss several examples.

\subsection{Families of contractions}
Let $T_1$,$\ldots $,$T_N$ be given contractions on the Hilbert space
$\cH $. For $\sigma =i_1^{k_1}\ldots i_n^{k_n}$,
$i_j\ne i_{j+1}$, $j=1,\ldots ,n-1$, $k_j\in \ZZ -\{0\}$, 
a reduced word in the free group $\FF _N$ on $N$ generators we define
the contraction
$$T_{\sigma }=T^{[k_1]}_{i_1}\ldots T^{[k_n]}_{i_n},$$
where 
$$T^{[k]}=\left\{\begin{array}{cc}
T^k & \quad k\geq 0;\\
T^{*-k} & \quad k<0.
\end{array}
\right.
$$
Then the kernel $K(\sigma ,\tau )=T_{\sigma ^{-1}\tau }$
is positive definite on $\FF _N$ (\cite{Boz}, \cite{Boc}).
It is also invariant, in the sense
that \eqref{invariance} holds for $\tau $, $\sigma $, $\sigma '$ in $\FF _N$.
Its restriction to $\FF ^+_N$, denoted $K^+$, 
is an invariant kernel on $\FF ^+_N$.
If we try to calculate the parameters of $K^+$ (with respect to the 
lexicographic order), we notice that their form is quite
complicated. However, it is not difficult to find the orthogonal polynomials. 
Thus, take $T_k=t_k$, $k=1,\ldots ,N$, where $t_k$ is a complex number in the
open unit disk $\DD $ and define $\phi :\cP _{2N}\rightarrow \CC $, 
\begin{equation}\label{dilat}
\phi (X_{i_1} \ldots X_{i_n})=t_{\epsilon (i_1)\ldots \epsilon (i_n)},
\end{equation}
where 
$$\epsilon (i)=\left\{\begin{array}{ccl}
i & \mbox {if} & 1\leq i\leq N; \\
 & & \\
(i-N)^{-1} & \mbox{if} & N<i\leq 2N 
\end{array}
\right. 
$$
(if $1\leq j\leq N$, then $j$ is viewed as an element of $\FF _N$ with 
inverse $j^{-1}$). Then 
$\phi (X^+_{\sigma }X_{\tau })=K^+(\sigma ,\tau )$ for 
$\sigma ,\tau \in \FF ^+_N$, so that $\phi $ is a  
q-positive functional with moment kernel $K^+$. Since $t_k\in \DD $, 
$k=1,\ldots N$, it follows that $\phi $ is a strictly q-positive
functional on $\cP _{2N}$ and since $G(\cP _{2N})=\FF ^+_N$, 
we can consider $\{\varphi _{\sigma }\}_{\sigma \in \FF ^+_N}$ 
the set of orthonormal polynomials associated with $\phi $.
\begin{theorem}\label{exemp}
The orthonormal polynomials of $\phi $ are: 
$$\varphi _k=
\frac{1}{d_{t_k}}(X_k-t_k), \quad k=1, \ldots N,$$
and for $\sigma \in \FF ^+_N$,
$$\varphi _{\sigma k}=X_{\sigma }\varphi _k, \quad k=1, \ldots N.$$
\end{theorem}
\begin{proof}
Let $\sigma ,\tau \in \FF ^+_N$, $\tau \prec \sigma $. Then 
$\sigma =\sigma 'k$ for some $\sigma '\in \FF ^+_N$ and some
$k=1,\ldots ,N$. Consequently, 
$$\begin{array}{rcl}
\phi (X^+_{\tau }\varphi _{\sigma })&=&
\phi (X^+_{\tau }X^+_{\sigma '}\varphi _k) \\
 & & \\
&=&\frac{1}{d_{t_k}}\left(\phi (X^+_{\tau }X^+_{\sigma '}X_k)-
t_k\phi (X^+_{\tau }X^+_{\sigma '})\right)\\
 & & \\
&=&\frac{1}{d_{t_k}}(t_{\tau ^{-1}\sigma 'k}-t_kt_{\tau ^{-1}\sigma '})=0.
\end{array}
$$
Also, 
$$\begin{array}{rcl}
\phi (\varphi ^+_{\sigma }\varphi _{\sigma })&=&
\phi (\varphi ^+_kX^+_{\sigma '}X_{\sigma '}\varphi _k) \\
 & & \\
&=&\frac{1}{d^2_{t_k}}\phi \left((X^+_k-\overline{t}_k)
X^+_{\sigma '}X_{\sigma '}(X_k-t_k)\right) \\
& & \\
&=&\frac{1}{d^2_{t_k}}(1-|t_k|^2-|t_k|^2+|t_k|^2)=1.
\end{array}
$$
These relations show that $\{\varphi _{\sigma }\}_{\sigma \in \FF ^+_N}$
is indeed the family of orthonormal polynomials associated with $\phi $.
\end{proof}

\noindent
Similar calculations will give an explicit formula of 
$\varphi _{\sigma }(X_1,\ldots ,X_N,l)$ for $l\geq 1$, at least for $|\sigma |$
large enough. 
Thus, we introduce the following notation: $r:\FF ^+_N\rightarrow \NN _0$ 
is the natural bijection between 
$\FF ^+_N$ and $\NN _0$, so that $r(\emptyset )=0$, $r(1)=1$, 
$\ldots $, $r(N)=N$, $r(11)=N+1$, $\ldots $; then 
for $l\geq 0$ and $\sigma \in \FF ^+$, $\sigma -l$ denotes the word in 
$\FF ^+_N$ that is $l$ steps ahead of $\sigma $ (so $\sigma -1$ is just
the predecessor of $\sigma $).
For $r(\sigma )>l$, the word $\sigma -l$ can be uniquely represented in 
the form $\sigma -l=q(\sigma )p(\sigma )$ for some 
$q(\sigma )\in \FF ^+_N$ and $p(\sigma )\in \{1,\ldots ,N\}$.
With this notation, we can obtain as in the proof of 
Theorem ~\ref{exemp} that for $l\geq 1$ and $r(\sigma )>l$, 
\begin{equation}\label{xeemp}   
\varphi ^l_{\sigma }=\frac{1}{d_{p(\sigma )}}\left(X_{\sigma }-
t_{p(\sigma )}X_{q(\sigma )}\right).
\end{equation}

As a consequence of \eqref{xeemp} and Theorem ~3.2 in \cite{BC}
we obtain $\gamma _{\sigma ,\tau }=0$ for $\sigma +r(\sigma )+l\prec 
\tau $. This gives more information about the parameters of $K^+$ but
still the remaining parameters look too complicated compared with 
the fact that $K^+$ is determined by just $N$ complex numbers.
A possibility to address this issue is to use parameters
associated to $K^+$ along a fixed chordal sequence, as in 
Theorem ~3.1 in \cite{BaC}.
More precisely, we use the following construction. 
In general we use the notation $G=(V, E)$ in order to denote an 
undirected graph with $V$ the set of vertices and $E$ the set of edges.
For $v,w\in V$, the notation $(v,w)$ denotes the edge of $G$ with 
endpoints $v$ and $w$.
Let $E^0_{\emptyset }=\emptyset $ and 
$G^0_{\emptyset }=(\FF ^+_N, \emptyset )$. 
For $\sigma \in \FF ^+_N-\{\emptyset \}$, $k\in \{1,\ldots ,N\}$,
and  $1\leq l\leq r(\sigma k)$, we define 
$$E^1_{\sigma k}=E^{r(\sigma k-1)}_{\sigma k-1}\cup 
\{(\sigma ,\sigma k)\}
$$
and for $l>1$, the set $E^l_{\sigma k}$ is obtained by adding one new
edge $(\tau ,\sigma k)$ to $E^{l-1}_{\sigma k}$, where $\tau \ne \sigma $ and 
$\tau \prec \sigma k$. Then define $G^l_{\sigma }=(\FF ^+_N,E^l_{\sigma k})$.
It is easily seen that $V_{\sigma }=\{\emptyset \preceq \tau \preceq \sigma \}$
is a maximal clique in $G^{r(\sigma )}_{\sigma }$, that is 
$(V_{\sigma },E^{r(\sigma )}_{\sigma })$ is the complete graph and 
$V_{\sigma }$ is maximal with this property. 
This implies that each $G^l_{\sigma }$ is a chordal graph and if we order the 
family $\{G^l_{\sigma }\}$ by lexicographic order on the pairs
$(\sigma ,l)$, $\sigma \FF ^+_N$, $1\leq l\leq r(\sigma )$, then 
$\{G^l_{\sigma }\}$ is a chordal sequence, according to the terminology 
in \cite{BaC}. By Theorem ~3.1 in \cite {BaC} (see also
Theorem ~ 7.2.7 in \cite{Co}), the kernel $K^+$ is uniquely determined by
a family $\{\gamma ^l_{\sigma } \mid \sigma \in \FF ^+_N-\{\emptyset\}, 
1\leq l\leq r(\sigma )\}$ of complex numbers with 
$|\gamma ^l_{\sigma }|<1$. We call these numbers the 
{\it parameters of $K^+$ along the chordal sequence} $\{G^l_{\sigma }\}$.

\begin{theorem}\label{mexep}  
The parameters of $K^+$ 
along the chordal sequence $\{G^l_{\sigma }\}$ are given by:
$\gamma ^1_{\sigma k}=t_k$ for $k=1,\ldots ,N$, $\sigma \in \FF ^+_N$ and 
$\gamma ^l_{\tau }=0$ for $l>1$ and $\tau \in \FF ^+_N-\{\emptyset\}$.
\end{theorem}
\begin{proof}
For $V\subset \FF ^+_N$ we denote by $K^+_V$ the restriction of $K^+$ to $V$,
that is, $K^+_V(\sigma ,\tau )=K^+(\sigma ,\tau )$
for $\sigma ,\tau \in V$. We claim that $K^+_{\{\emptyset \preceq 
\tau \preceq \sigma k\}}$ is the Markov product of the kernels
$K^+_{\{\emptyset \preceq \tau \preceq \sigma k-1\}}$ and 
$K^+_{\{\sigma ,\sigma k\}}$. Indeed, we have
$\{\emptyset \preceq \tau \preceq \sigma k-1\} \cap \{\sigma ,\sigma k\}
=\{\sigma \}$
and 
$$\begin{array}{rcl}
K^+(\tau ,\sigma k)&=&t_{\tau ^{-1}\sigma k}=t_{\tau ^{-1}\sigma }t_k
=t_{\tau ^{-1}\sigma }t_{\sigma ^{-1}\sigma k} \\
& & \\
&=& K^+(\tau, \sigma )K^+(\sigma ,\sigma k).
\end{array}
$$
Now an application of Theorem ~2.3 concludes the proof.
\end{proof}

We could deal now with orthogonal polynomials along a chordal sequence
such as the one above. However, we do not pursue this here, more
details can be found in \cite{Banks}.

\subsection{Moment kernels on $\cA ^N_O$}
We can see that the functional $\phi $ given by 
\eqref{dilat} induces a functional $\tilde \phi $ on $\cA ^N_O$
such that $\tilde \phi \circ \pi _{\cA ^N_O}=\phi $.
This suggests that the invariant kernels are related to $\cA ^N_O$ and the 
following result explains this connection. We use the notation introduced in 
Subsection~2.1.

\begin{theorem}\label{motiv}
$K=K_{\phi }$
for some linear functional on $\cR (\cA ^N_O)$ if and only if
$K$ is an invariant kernel.
\end{theorem}
\begin{proof}
Let $K=K_{\phi }$
for some linear functional $\phi $ on  $\cR (\cA ^N_O)$ and let
$\tau, \sigma ,
\sigma '$ be words in the index set of $\cA ^N_O$, which is $\FF ^+_N$.
Then
$$\begin{array}{rcl}
K(\tau \sigma ,\tau \sigma ') &= 
& \phi (X^+_{\tau \sigma }X_{\tau \sigma '}) \\
 & & \\
& =&\phi (X^+_{\sigma }X^+_{\tau }X_{\tau }X_{\sigma '}.
\end{array}
$$
Since $X^+_{\tau }X_{\tau }=1$ in $\cR (\cA ^N_O)$, we deduce that
$$K(\tau \sigma ,\tau \sigma ')=K(\sigma ,\sigma ').$$

Conversely, let $K$ be an invariant kernel. Any element of $\cR(\cA ^N_O)$ 
is a linear combination of monomials $X_{i_1}\ldots X_{i_n}$,
$i_1$, $\ldots $, $i_n\in \{1,\ldots ,2N\}$, with the 
property that there is no pair $(i_k,i_{k+1})$ with $i_k>N$ 
and $i_k-i_{k+1}=N$. We define $\phi $ on monomials as above which can 
be written in the form $X^+_{\sigma }X_{\tau }$ by the formula
$$\phi (X^+_{\sigma }X_{\tau })=K(\sigma ,\tau ),$$
and arbitrarly on the other monomials in $\cR (\cA ^N_O)$. 
The invarince of $K$ insures that $\phi $ is well-defined. Then we
extend $\phi $ by linearity to the whole $\cR (\cA ^N_O)$ and 
clearly $K=K_{\phi }$. 
\end{proof}

This result explains that the study of orthogonal polynomials on 
$\cR (\cA ^N_O)$ reduces to the study of invariant kernels.

\subsection{Free products}
Since $\cP _N$ is a free product of $N$ copies of $\cP _1$, 
it is quite natural to look at free products of q-positive
functionals. 
Let $\cR (\cA _1)$, $\cR (\cA _2)$ be two algebras with sets of defining
relations $\cA _1$, respectively, $\cA _2$. It is convenient to view
$\cR (\cA _1)$ as a quotient of $\cP_{2N}$ in the variables $X_1$,
$\ldots $, $X_{2N}$ and $\cR (\cA _2)$ as a quotient of 
$\cP_{2M}$ in the variables $Y_1$,
$\ldots $, $Y_{2M}$. According to the notation in Subsection ~2.1, 
let $F_{\alpha }=\pi _{\cA _1}(X_{\alpha })$, 
$\alpha \in G(\cA _1)$, and 
$G_{\beta }=\pi _{\cA _2}(Y_{\beta })$, 
$\beta \in G(\cA _2)$. Each of $G(\cA _1)$ and 
$G(\cA _2)$ contains words of length $1$, otherwise the situation is 
degenerate, in the sense that $\cR (\cA _1)=\cR (\cA _2)=\CC $. In order to
simplify the notation, but without loss of generality, we can assume that 
$G(\cA _1)$ contains all of $1$, $\ldots $, $N$ and  
$G(\cA _2)$ contains all of $1$, $\ldots $, $M$. In this way, 
$\cR (\cA _1)$ is the set of polynomials in the variables 
$F_1$, $\ldots $, $F_N$, $F^+_1$, $\ldots $, $F^+_N$
(satisfying the defining relations in $\cA _1$), and similarly, 
$\cR (\cA _2)$ is the set of polynomials in the variables 
$G_1$, $\ldots $, $G_M$, $G^+_1$, $\ldots $, $G^+_M$
(satisfying the defining relations in $\cA _2$). 
Let $\cR ^0(\cA _i)$, $i=1,2$, denote the set of polynomials in 
$\cR (\cA _i)$ without constant term. Then 
$$\cR _1(\cA _1)\star \cR _1(\cA _2)=
\CC \oplus \left(\oplus _{n\geq 1}
\oplus _{i_1\ne i_2,\ldots ,i_{n-1}\ne i_n}\cR ^0(\cA _{i_1})
\otimes \ldots \otimes \cR ^0(\cA _{i_n})\right),
$$
and we notice that $\cR _1(\cA _1)\star \cR _1(\cA _2)$ is isomorphic to 
$\cR (\cA _1+\cA _2)$, where $\cA _1+\cA _2$ is the disjoint union of 
$\cA _1 $ and $\cA _2$ (due to our convention to view $\cA _1$ as a subset 
of $\cP _{2N} $ in the variables  
$X_1$,
$\ldots $, $X_{2N}$ and $\cA _2$ as a subset of 
$\cP_{2M}$ in the variables $Y_1$,
$\ldots $, $Y_{2M}$, the sets $\cA _1 $ and $\cA _2$ 
are automatically disjoint).

Now let $\phi _1 $ be a q-positive functional on $\cR (\cA _1)$ and 
$\phi _2 $ be a q-positive functional on $\cR (\cA _2)$. Their free
product $\phi =\phi _1\star \phi _2$ on $\cR _1(\cA _1)\star \cR _1(\cA _2)$
is defined by $\phi (1)=1$ and 
$\phi (P_{i_1}\ldots P_{i_n})=\phi _{i_1}(P_{i_1})\ldots 
\phi _{i_n}(P_{i_n})$ for $n\geq 1$, 
$i_1\ne i_2$, $\ldots $, $i_{n-1}\ne i_n$, $P_{i_k}\in \cR ^0(\cA _{i_k})$,
and $i_k\in \{1,2\}$ for $k=1,\ldots ,n$. 
The map $\phi $ given by \eqref{dilat} is an example of a free product
of $N$ q-positive functionals.

Since in general a q-positive functional is not positive, the main result in 
\cite{Boz} cannot be applied in order to conclude
that the free product of two q-positive functionals is q-positive, 
however this 
follows from the more general result in \cite{Boc}.
Using Theorem ~\ref{motiv} it follows that for 
strictly q-positive functionals $\phi _1$ and $\phi _2$
on $\cR (\cA ^N_0)$ and, respectively, $\cR (\cA ^M_0)$, 
the kernel $K_{\phi _1\star \phi _2 }$ is a strictly positive definite 
invariant
kernel.
This construction produces a relatively large class of 
positive definite and strictly positive definite invariant kernels.

\section{The displacement structure of invariant kernels}

From the examples in the previous section we see that it is difficult to 
explore the additional symmetry of a positive definite invariant kernel 
in terms of its parameters or of its orthogonal polynomials. In 
particular, the invariance is not encoded efficiently into the 
generators of the displacement equation   
\eqref{displ}. In this section we consider a different 
displacement structure of an invariant kernel. In order to avoid
notational complications, we can assume that the positive definite
invariant kernel $K$ is scalar-valued, and also that 
$K(\sigma ,\sigma )=1$ for all $\sigma \in \FF ^+_N$. For each $n\geq 0$ we 
introduce the following elements: the $\sum _{k=0}^nN^k\times  
\sum _{k=0}^nN^k$ matrix $F_{k,n}$, $k=1, \ldots ,N$, whose action 
on the Hilbert space $\cF _n$ of sequences 
$\{h_{\sigma }\}_{|\sigma |\leq n}$ (with Euclidean norm) is given by
$$F_{k,n}(\{h_{\sigma }\}_{|\sigma |\leq n})=
\{g_{\sigma }\}_{|\sigma |\leq n},$$
where 
$$g_{\tau }=\left\{
\begin{array}{rl}
h_{\sigma } & \mbox{if $\tau =k\sigma $ } \\
 & \\
 0 & \mbox{otherwise}.
\end{array}
\right.
$$
Also, let 
$R_n=\left[K(\sigma ,\tau )\right]_{|\sigma |,|\tau |\leq n}$
and define $Q_n=\left[Q_n(\sigma ,\tau )\right]_{|\sigma |,|\tau |\leq n}$,
where $Q_n(\sigma ,\tau )=0$ if  
$\sigma =\alpha \sigma '$, $\tau =\alpha \tau '$
for some $\alpha \in \FF ^+_N-\{\emptyset \}$, $\sigma ',\tau '\in \FF ^+_N$, 
and  
otherwise $Q_n(\sigma ,\tau )=K(\sigma ,\tau )$. The next result
shows that the left hand side of the relation \eqref{patagonia}
sifts out all the redundancy in $K$ caused by its invariance.

\begin{lemma}\label{final}
For each $n\geq 0$ the matrix $R_n$ satisfies the displacement 
equation
\begin{equation}\label{patagonia}
R_n-\sum _{k=1}^NF_{k,n}R_nF^*_{k,n}=Q_n.
\end{equation}
\end{lemma}
\begin{proof}
Let $\{e_{\sigma }\}_{|\sigma |\leq n}$ be the standard basis 
of the Hilbert space $\cF _n$. Then
$$\begin{array}{rl}
\langle (R_n-\sum _{k=1}^NF_{k,n}R_nF^*_{k,n})e_{\sigma },e_{\tau}\rangle & \\
 & \\
& \!\!\!\!
\!\!\!\!
\!\!\!\!
\!\!\!\!\!\!\!\!
\!\!\!\!
\!\!\!\!
\!\!\!\!\!\!\!\!
\!\!\!\!
\!\!\!\!
\!\!\!\!
=\langle R_ne_{\sigma },e_{\tau}\rangle -
\sum _{k=1}^N\langle F_{k,n}R_nF^*_{k,n}e_{\sigma },e_{\tau}\rangle  \\
  & \\
& \!\!\!\!
\!\!\!\!
\!\!\!\!
\!\!\!\!\!\!\!\!
\!\!\!\!
\!\!\!\!
\!\!\!\!\!\!\!\!
\!\!\!\!
\!\!\!\!
\!\!\!\!
=\langle R_ne_{\sigma },e_{\tau}\rangle -
\sum _{k=1}^N\langle F_{k,n}R_nF^*_{k,n}e_{\sigma },e_{\tau}\rangle .
\end{array}
$$
If there is no  $\alpha \in \FF ^+_N-\{\emptyset \}$ such that
$\sigma =\alpha \sigma '$ and $\tau =\alpha \tau '$,
then the first letter of $\sigma $ is going to be different 
from the first letter of $\tau $, which implies that 
$\sum _{k=1}^N\langle F_{k,n}R_nF^*_{k,n}e_{\sigma },e_{\tau}\rangle =0$
and 
$$\begin{array}{rl}
\langle (R_n-\sum _{k=1}^NF_{k,n}R_nF^*_{k,n})e_{\sigma },e_{\tau}\rangle & \\
 &  \\
 &\!\!\!\!
\!\!\!\!
\!\!\!\!
\!\!\!\!\!\!\!\!
\!\!\!\!
\!\!\!\!
\!\!\!\!\!\!\!\!
\!\!\!\!
\!\!\!\!
\!\!\!\!
=\langle R_ne_{\sigma },e_{\tau}\rangle =K(\tau ,\sigma ) \\
 & \\ 
 &\!\!\!\!
\!\!\!\!
\!\!\!\!
\!\!\!\!\!\!\!\!
\!\!\!\!
\!\!\!\!
\!\!\!\!\!\!\!\!
\!\!\!\!
\!\!\!\!
\!\!\!\!
=Q_n(\tau ,\sigma )=\langle Q_ne_{\sigma },e_{\tau}\rangle .
\end{array}
$$
If there is $\alpha \in \FF ^+_N-\{\emptyset \}$ such that
$\sigma =\alpha \sigma '$ or $\tau =\alpha \tau '$, this implies that there
is $p\in \{1, \ldots ,N\}$ such that $\sigma =p\sigma '$, 
$\tau =p\tau '$, and then 
$$\begin{array}{rl}
\langle (R_n-\sum _{k=1}^NF_{k,n}R_nF^*_{k,n})e_{\sigma },e_{\tau }\rangle & \\
  & \\
& \!\!\!\!
\!\!\!\!
\!\!\!\!
\!\!\!\!\!\!\!\!
\!\!\!\!
\!\!\!\!
\!\!\!\!\!\!\!\!
\!\!\!\!
\!\!\!\!
\!\!\!\!
=K(p\tau ',p\sigma ')-
\sum _{k=1}^N\langle F_{k,n}R_nF^*_{k,n}e_{p\sigma '},e_{p\tau '}\rangle   \\
 &  \\
& \!\!\!\!
\!\!\!\!
\!\!\!\!
\!\!\!\!\!\!\!\!
\!\!\!\!
\!\!\!\!
\!\!\!\!\!\!\!\!
\!\!\!\!
\!\!\!\!
\!\!\!\!
=K(\tau ',\sigma ')-\langle R_ne_{\sigma '},e_{\tau '}\rangle =
0=Q_n(\tau ,\sigma ).
\end{array} 
$$
In conclusion, we obtained \eqref{patagonia}.
\end{proof}

We now try to factorize $Q_n$ in the form $G_nJ_nG^*_n$
for some symmetry $J_n$ ($J_n=J_n^*=J^{-1}_n$), but of course, 
$J_n$ is no longer
$\left[
\begin{array}{cc}
1 & 0 \\
0 & -1 
\end{array}
\right]
$. 
In order to obtain a result suitable for the displacement structure
theory, $J_n$ should be the same for all $Q_n$
(that is, for all invariant kernels $K$).

\begin{lemma}\label{suspanse}
$(a)$ Let $A=\left[A_{i,j}\right]_{i,j=1}^p$ be a selfadjoint block-matrix
with $A_{k,k}=0$ for all $k=1,\ldots ,p$. Then 
\begin{equation}\label{factorizare}
A=B\fI _{2p-2}B^*,
\end{equation}
where 
\begin{equation}\label{bi}
B=\left[\begin{array}{ccccccccc}
0 & 0 & \ldots & 0 & 0 & 0 & \ldots & 0 & I \\
A_{2,1} & 0 & & & 0 & 0 & & I & 0 \\
A_{3,1} & A_{3,2} & & & 0 & 0 & & 0 & 0 \\
\vdots & \vdots & \ddots & & \vdots & \vdots & & \vdots & \vdots \\
A_{p-1,1} & A_{p-1,2} & \ldots & 0 & I & 0 & \ldots & 0 & 0 \\
A_{p,1} & A_{p,2} & \ldots & A_{p,p-1} & 0 & 0 & \ldots & 0 & 0
\end{array}
\right]
\end{equation}
is a $p\times (2p-2)$ block matrix and $\fI _k$ is a $k\times k$
block-matrix,
\begin{equation}\label{ji}
\fI _k=\left[
\begin{array}{ccccc}
0 & 0 & \ldots & 0 & I \\
0 & 0 & \ldots & I & 0 \\
\vdots & \vdots & & \vdots & \vdots \\
0 & I & \ldots & 0 & 0 \\
I & 0 & \ldots & 0 & 0 
\end{array}
\right].
\end{equation}

$(b)$ Assume all $A_{i,j}$ are complex numbers. Then 
$2p-2$ is the minimal dimension of a symmetry $J$ 
with the property that 
for any  selfadjoint matrix $A$ with zero diagonal,
the relation 
\eqref{factorizare}
holds for some $B$.
\end{lemma}
\begin{proof}
$(a)$ The formula \eqref{factorizare} is easily verified by 
direct computations 
that can be omitted.

$(b)$ Since $A_{k,k}=0$ for $k=1, \ldots ,p$, $A$ cannot be 
positive or negative (excepting for the trivial case $A_{i,j}=0$ for all 
$i,j$). Therefore $A$ can have (generically) at most 
$p-1$ positive eigenvalues or at most $p-1$ negative eigenvalues.
This implies that the symmetry $J$ satisfying \eqref{factorizare} 
for all selfadjoint $A$ with zero diagonal must have at least as many
positive eigenvalue and, respectively,  negative eigenvalues, 
which gives a total
of at least $2p-2$ eigenvalues. The construction of $(a)$ realizes this value, 
so $2p-2$ is the minimal dimension of a symmetry $J$ satisfying 
\eqref{factorizare} for any selfadjoint matrix with zero diagonal.
\end{proof} 

\begin{theorem}\label{realfinal}
For each $n\geq 0$ the matrix $R_n$ satisfies the displacement 
equation
\begin{equation}\label{realpatagonia}
R_n-\sum _{k=1}^NF_{k,n}R_nF^*_{k,n}=G_nJ_nG^*_n,
\end{equation}
where $J_n$ is a symmetry of dimension 
$2+(2N-2)\sum _{k=0}^{n-1}N^k$.
\end{theorem}
\begin{proof}
By Lemma ~\ref{final}, the matrix $R_n$ satisfies the displacement 
equation 
$$R_n-\sum _{k=1}^NF_{k,n}R_nF^*_{k,n}=Q_n.
$$
From the definition of $Q_n$, we deduce that
$$Q_n=\left[\begin{array}{cc}
1 & S_n \\
S^*_n & 0 
\end{array}
\right]+ \left[\begin{array}{cc}
0 & 0 \\
0 & L_n 
\end{array}
\right],$$
wher $S_n=\left[K(\emptyset ,\sigma )\right]_{|\sigma |\leq n}$, 
$L_n(\sigma ,\tau )=0$ if $\sigma =\alpha \sigma '$, 
$\tau =\alpha \tau '$ for some $\alpha \in \FF ^+_N-\{\emptyset \}$, 
$\sigma ',\tau '\in \FF ^+_N$, and otherwise, 
$L_n(\sigma ,\tau )=K(\sigma ,\tau )$ (note that $L_n(\sigma ,\tau )$
is defined only for $\sigma ,\tau \in  \FF ^+_N-\{\emptyset \}$).
Since we have the factorization
$$\left[\begin{array}{cc}
1 & S_n \\
S^*_n & 0 
\end{array}
\right]=\left[\begin{array}{cc}
1 & 0 \\
S^*_n & S^*_n   
\end{array}
\right]
\left[\begin{array}{cc}
1 & 0 \\
0 & -1
\end{array}
\right]
\left[\begin{array}{cc}
1 & S_n \\
0 & S_n 
\end{array}
\right], 
$$
we need only to show that $L_n$ has a factorization 
of the form $L_n=G'_nJ'_nG'^*_n$ with a symmetry
$J'_n$ of dimension $(2N-2)\sum _{k=0}^{n-1}N^k$. Then 
$$Q_n=\left[\begin{array}{ccc}
1 & 0 & 0 \\
S^*_n & S^*_n & G'_n
\end{array}
\right]
\left[\begin{array}{ccc}
1 & 0 & 0 \\
0 & -1 & 0 \\
0 & 0 & J'_n 
\end{array}
\right]
\left[\begin{array}{cc}
1 & S_n \\
0 & S_n \\
0 & G'^*_n
 \end{array}
\right]
$$
will be the required factorization of $Q_n$.

Now, for $k=1, \ldots ,N$ we define 
$\cA _k=\{k\tau \mid |\tau |\leq n-1 \}.$
Then $\{\cA _k\}_{k=1}^N$ is a partition of the set
$\cW _n=\{\sigma \in \FF ^+_N \mid |\sigma |\leq n\}$. We reorder 
the elements of $\cW _n$ such that $\sigma <\tau $ if 
$\sigma \in \cA _k$, $\tau \in \cA _j$, $k<j$.
Then $Q_n=\left[A_{i,j}\right]_{i,j=1}^N$ with $A_{k,k}=0$ for all 
$k=1,\ldots ,N$. By Lemma \ref{suspanse},
$$Q_n=B_n\fI _{2N-2}B^*_n,$$
where $B_n$ is given by \eqref{bi}. We can define $J'_N=\fI _{2N-2}$, 
so that we obtain a factorization of $Q_n$  with the required 
dimension of the symmetry $J_n$. 
\end{proof}

\noindent
The symmetry $J_n$ in Theorem \ref{realfinal} is unitarily 
equivalent to the symmetry
$\left[\begin{array}{cc}
I_{p_n} & 0 \\
0 & -I_{p_n}
\end{array}
\right]$, where 
$p_n=1+(N-1)\sum _{k=0}^{n-1}N^k$, so that we can rewrite 
equation \eqref{realpatagonia}
in the more familiar form 
$$R_n-\sum _{k=1}^NF_{k,n}R_nF^*_{k,n}=G_n
\left[\begin{array}{cc}
I_{p_n} & 0 \\
0 & -I_{p_n}
\end{array}
\right]G^*_n,
$$
for some new $G_n$, and the established results of the displacement structure
theory can be used in order to explore the structure of $R_n$. In 
particular, we obtain a Schur type algorithm that better encodes the 
invariance of $K$. Still, we have to note the fact that due to the 
complexity of $K$, the number $p_n$ depends on $n$. Some more details in 
this direction can be found in \cite{Banks}.

We conclude by noticing that the moment kernel of a q-positive functional
on  $\cR (\cA ^N_{CT})$ is characterized by the property that 
$L_n=0$ for all $n\geq 0$. Thus, 
\eqref{realpatagonia} appears as an extension of the displacement 
equation for $\cR (\cA ^N_{CT})$ obtained in \cite{CJ}.

\end{document}